\documentclass[12pt,final]{amsart}

\usepackage[margin=1in]{geometry}

\usepackage{amssymb}

\newtheorem{lemma}{Lemma}[section]

\newtheorem{theorem}[lemma]{Theorem}

\newtheorem{corollary}[lemma]{Corollary}

\theoremstyle{definition}

\theoremstyle{remark}
\newtheorem{remark}[lemma]{Remark}

\newcommand{\B}{\ensuremath{{\mathbb B}}}
\newcommand{\N}{\ensuremath{{\mathbb N}}}

\newcommand{\R}{\ensuremath{\mathbb R}}
\newcommand{\C}{\ensuremath{\mathbb C}}

\newcommand{\K}{\ensuremath{\mathcal{K}}}

\author{Tatyana Barron}
\address{T. Barron, Department of Mathematics, 
University of Western Ontario, 
London, Ontario N6A 5B7, Canada }
\email{tatyana.barron@uwo.ca}

\thanks{Research is supported in part 
by the Natural Sciences and Engineering Research Council of Canada}

\title[Closed geodesics and pluricanonical sections]{Closed geodesics and pluricanonical sections on ball quotients}

\begin{document}
\sloppy

\maketitle

\noindent {\bf Abstract.}  We obtain asymptotics of sequences of 
the holomorphic sections of the pluricanonical bundles on ball quotients associated 
to closed geodesics. A nonvanishing result follows.

\

\noindent {\bf Keywords:} canonical bundle, holomorphic section, asymptotics, complex hyperbolic space, geodesic. 

\

\section{Introduction}

Let $n\in\N$ and 
let $\Gamma$ be a cocompact discrete subgroup of $SU(n,1)$ that acts freely on the unit ball $\B^n$ in $\C^n$. 
We consider the ball with the invariant Bergman metric. 
Let $k\ge 2$ be 
an integer. Suppose $\gamma\in\Gamma$, $\gamma\ne id$, is such that all of its eigenvalues are real, 
and moreover, the endpoints of the $\gamma$-invariant geodesic $\tilde{C}$ on $\partial \B^n$ are also real. We describe 
a construction of a holomorphic section of the $k$-th tensor power of the canonical bundle on $\B^n/\Gamma$   
that is associated to the closed geodesic $C=\tilde{C}/\langle \gamma\rangle$ in $\B^n/\Gamma$, we determine the $k\to\infty$ asymptotics 
of the norm of this section, and we conclude that the norm is not zero for large $k$ (Theorem \ref{mainth}, Corollary \ref{nonvcor}). 

Associating sections of vector bundles to submanifolds can be done in different ways, for a variety of purposes.  
See, for this kind of ideas applied to ball quotients and geodesic cycles 
\cite{kudla:79, tong:83}, and also \cite{borthwick:95, foth:01, foth:04, katok:85, katok:87}. 
Another frequently used approach is associating a section of a line bundle on a compact K\"ahler manifold to a compact Lagrangian 
submanifold. See, in particular, references \cite{borthwick:95, burns:10, deber:06, jeffrey:92, paoletti:08} 
and the papers \cite{foth:02, foth:08} by T. Foth (T. Barron). 
In our setting described above, for $n>1$, $C$ is an isotropic submanifold of $\B^n/\Gamma$, but not a Lagrangian submanifold. 
In subsection \ref{sec:addrem} we explain more carefully how our main theorem fits into the general context.

\section{Preliminaries}

In this section we collect standard facts, for which the references are 
\cite{goldman:99,krantz:01,rudin:80,selberg:57,zhu:05}. 

Let $z_1$,..., $z_n$ be the complex coordinates on $\C^n$, $z_j=x_j+iy_j$ ($x_j,y_j\in\R$, $j=1,...,n$). 
The Hermitian symmetric space $SU(n,1)/S(U(n)\times U(1))$ is isomorphic to the unit ball  
$\B^n=\{ (z_1,...,z_n)\in\C^n|  \ |z_1|^2+...+|z_n|^2<1\}$.   
The group $SU(n,1)$ acts on $\B^n$  by fractional-linear transformations. 
We will use the same letter to denote a matrix from $SU(n,1)$ and the corresponding automorphism of the ball.   
For $A =(a_{jk})\in SU(n,1)$ the corresponding biholomorphism 
$\B^n \to \B^n$ is the mapping 
$$
z=(z_1,...,z_n)\mapsto \Bigl ( \frac{a_{11}z_1+...+a_{1n}z_n+a_{1,n+1}}{a_{n+1,1}z_1
+...+a_{n+1,n}z_n+a_{n+1,n+1}} ,...,  
\frac{a_{n1}z_1+...+a_{nn}z_n+a_{n,n+1}}{a_{n+1,1}z_1+
...+a_{n+1,n}z_n+a_{n+1,n+1}} \Bigr ).
$$
The complex Jacobian of the transformation $A$ at $z\in \B^n$ is 
\begin{equation}
\label{jacobidet}
J(A,z)=\frac{1}{(a_{n+1,1}z_1+...+a_{n+1,n}z_n+a_{n+1,n+1})^{n+1}}.
\end{equation}
For $\gamma_1,\gamma_2\in SU(n,1)$ and $z\in\B^n$ 
\begin{equation}
\label{cocyclecond}
J(\gamma_1\gamma_2,z)=J(\gamma_1,\gamma_2 z)J(\gamma_2,z).
\end{equation}
The group $SU(n,1)$ also acts on $\partial\B^n=\{ (z_1,...,z_n)\in\C^n|  \ |z_1|^2+...+|z_n|^2=1\}$, 
by fractional-linear transformations. 

For two vectors $u,v$ in $\C^{n+1}$ denote 
$$
\langle \langle u,v\rangle \rangle =u_1\bar v_1+...+u_n\bar v_n-u_{n+1}\bar v_{n+1}=u^T\sigma \bar v,
$$
where
$$
\sigma=\begin{pmatrix}
1_{n} & 0 \\
0 & -1
\end{pmatrix}.
$$
For $z,w\in \B^n$ denote 
$$
\langle z,w\rangle=z_1\bar{w}_1+...+z_n\bar{w}_n-1
$$
and denote by $\rho (z,w)$ the distance between $z$ and $w$ with respect 
to the complex hyperbolic metric. Note that 
$$
\cosh ^2\frac{\rho(z,w)}{2}=\frac{\langle z,w\rangle \langle w,z\rangle}{\langle z,z\rangle \langle w,w\rangle}. 
$$ 
Let $g$ be a loxodromic element of $SU(n,1)$ (i.e. $g$ has exactly two fixed points in $\B^n\cup\partial \B^n$, 
these two fixed points will necessarily be on $\partial \B^n$). 
We will call $g$ {\it hyperbolic} if all of its eigenvalues are real. Denote the $n+1$ eigenvalues of $g$ by 
 $\alpha_j$, $1\le j\le n+1$. Let $\alpha_n$ and $\alpha_{n+1}$ be the eigenvalues for the eigenvectors of $g$ 
that correspond to the fixed points $X,Y\in \partial \B^n$, respectively. Then $\alpha_{n+1}=\dfrac{1}{\alpha_n}$, and,  without loss of generality,   
$|\alpha_n|>1$.  Each of the eigenvalues $\alpha_1$,...,$\alpha_{n-1}$ 
is either $1$ or $-1$. If both $1$ and $-1$ occur among the eigenvalues,  then the ordering is assumed to be so that 
$\alpha_1=...=\alpha_{2l}=-1$ and $\alpha_{2l+1}=...=\alpha_{n-1}=1$ for $l$ between $1$ and $n-1$.

So, the vectors  $\begin{pmatrix} X \\ 1\end{pmatrix}$ and $\begin{pmatrix} Y \\ 1\end{pmatrix}$ 
are eigenvectors of $g$ for $\alpha_n$ and $\alpha_{n+1}$ respectively. 
Let $v_1$,...,$v_{n-1}$ be the eigenvectors of $g$ corresponding to $\alpha_1$,...,$\alpha_{n-1}$, such that  
$v_j$ which are the eigenvectors for the eigenvalue $-1$  
form an orthonormal basis 
(with respect to $\langle \langle .,.\rangle \rangle$) in their linear span, 
and  $v_j$ which are the eigenvectors for the eigenvalue $1$ 
form an orthonormal basis 
(with respect to $\langle \langle .,.\rangle \rangle$) in their linear span. 
For each $j\in \{ 1,....,n-1\}$ the vector $v_j$ is orthogonal to $\begin{pmatrix} X \\ 1\end{pmatrix}$ and is orthogonal to 
$\begin{pmatrix} Y \\ 1\end{pmatrix}$
(with respect to $\langle \langle .,.\rangle \rangle$). If $v_j$ is an eigenvector for eigenvalue $1$ and $v_r$ 
is an eigenvalue for eigenvalue $-1$, then $\langle \langle v_j,v_r\rangle \rangle=0$. We also have: 
$\langle X,X\rangle=0$,  $\langle Y,Y\rangle=0$, $\langle X,Y\rangle\ne 0$. 
The matrix 
$$
A_g=\begin{pmatrix} v_1 & ... & v_{n-1} & \dfrac{1}{\langle X,Y\rangle}\begin{pmatrix} X \\ 1\end{pmatrix}+\dfrac{1}{2}\begin{pmatrix} Y \\ 1\end{pmatrix} & 
\dfrac{1}{\langle X,Y\rangle}\begin{pmatrix} X \\ 1\end{pmatrix}-\dfrac{1}{2}\begin{pmatrix} Y \\ 1\end{pmatrix} \end{pmatrix}
$$
is in $SU(n,1)$. The corresponding automorphism of $\B^n\cup\partial \B^n$ maps $(0,...,0,1)$ to $X$, 
maps $(0,...,0,-1)$ to $Y$. Let  
$\tilde{C}_0$ be the geodesic in $\B^n$ that connects $(0,...0,-1)$ and $(0,...,0,1)$:  
\begin{equation}
\label{tC0def}
\tilde{C}_0=\{ (z_1,...,z_n)\in \B^n| z_1=...=z_{n-1}=0, z_n=x_n+iy_n, y_n=0, -1<x_n<1\} .
\end{equation}
The transformation $A_g$ maps $\tilde{C}_0$ to the geodesic connecting $X$ and $Y$. Also $A_g^{-1}=\sigma \bar A_g^T \sigma$ and 
$$
A_g^{-1}gA_g=\begin{pmatrix} I_{g} & 0 & 0 \\ 
0 & \dfrac{1}{2\alpha_n}+\dfrac{\alpha_n}{2} & -\dfrac{1}{2\alpha_n}+\dfrac{\alpha_n}{2} \\ 
0 & -\dfrac{1}{2\alpha_n}+\dfrac{\alpha_n}{2} & \dfrac{1}{2\alpha_n}+\dfrac{\alpha_n}{2} 
\end{pmatrix}
$$
where $I_{g}=\begin{pmatrix} \alpha_1 & & \\ & ... & \\ & & \alpha_{n-1}\end{pmatrix}$.

The Euclidean volume form on $\B^n$ is 
$$
dV_e=dx_1\wedge dy_1\wedge...\wedge  dx_n\wedge dy_n=
\Bigl( \frac{i}{2}\Bigr )^n dz_1\wedge d\bar{z}_1\wedge ...\wedge dz_n\wedge d\bar{z}_n.
$$   
The Bergman kernel 
for $\B^n$ is 
$$
\K (z,w)=\frac{n!}{\pi^n}\frac{1}{(-\langle z,w\rangle)^{n+1}}.
$$ 
It has the reproducing property: 
$$
f(z)=\int\limits_D f(w) \K (z,w)dV_e (w),
$$
$z\in \B^n$, for all functions $f$ that are holomorphic on $\B^n$ and such that $\int\limits_{\B^n} |f(z)|^2dV_e(z)<\infty$. Also 
$\K (z,w)=\overline{\K (w,z)}$ for $z,w\in \B^n$, and 
\begin{equation}
\label{bktransf}
J(A,z) \overline{J(A,w)}\K (A z,A w)=\K (z,w)
\end{equation}
for $z,w\in \B^n$, $A \in SU(n,1)$. 
The K\"ahler form $i\partial\bar{\partial}\log K(z,z)$ and the volume form 
$$
dV(z)=\K (z,z) dV_e(z)
$$ 
on $\B^n$ are $SU(n,1)$-invariant.

Let $k$ be a positive integer. 
The reproducing kernel for the Hilbert space 
of holomorphic functions on $\B^n$ satisfying $\int\limits_{\B^n} |f(z)|^2\K (z,z)^{-k}dV(z)<\infty$ 
is $c(\B^n,k)\K (z,w)^k$, where $c(\B^n,k)=\binom{(n+1)(k-1)+n}{n}$. 
The reproducing property is, for any such function $f$: 
\begin{equation}
\label{wreprod}
f(z)=c(\B^n,k)\int\limits_{\B^n} f(w) \K (z,w)^k\K (w,w)^{-k}dV (w),
\end{equation}
$z\in \B^n$. 
Using the Stirling formula (see e.g. \cite{debru:61}), we get:  
\begin{equation}
\label{ballconstk}
c(\B^n ,k)\sim \frac{(n+1)^n k^n}{n!}\Bigl ( 1+O(\frac{1}{k})\Bigr )
\end{equation}
as $k\to\infty$.

\section{Asymptotics}

\subsection{The setting and the main result}
\label{sec:setting} 
Let $\Gamma$ be a discrete subgroup of $SU(n,1)$ such that $M=\Gamma\backslash \B^n$ is smooth and compact.   
Then all elements of $\Gamma$ are loxodromic. 

Assume that $\Gamma$ contains a hyperbolic element $\gamma$. Let $\tilde{C}$ be the geodesic in $\B^n$ invariant under $\gamma$. 
Assume that the endpoints $X$ and $Y$  of this geodesic in $\partial \B^n$ are real.  
The vectors $\begin{pmatrix} X \\ 1\end{pmatrix}$ and $\begin{pmatrix} Y \\ 1\end{pmatrix}$ 
are eigenvectors of $\gamma$ with eigenvalues $\lambda$ and $1/\lambda$ respectively, 
for some $\lambda\in\R$ ($\lambda\ne 0,\pm 1$). Assume that $\gamma$ is  not a power of any other element in $\Gamma$. 
Also assume that $\lambda^2 >1$ (if $\lambda^2<1$, then we can replace $\gamma$ by $\gamma^{-1}$). 
The automorphism of $\B^n\cup \partial \B^n$ defined by the matrix 
\begin{equation}
\label{Agammamatrix}
A_{\gamma}=
\begin{pmatrix} v_1 & ... & v_{n-1} & \dfrac{1}{\langle X,Y\rangle}\begin{pmatrix} X \\ 1\end{pmatrix}+\dfrac{1}{2}\begin{pmatrix} Y \\ 1\end{pmatrix} & 
\dfrac{1}{\langle X,Y\rangle}\begin{pmatrix} X \\ 1\end{pmatrix}-\dfrac{1}{2}\begin{pmatrix} Y \\ 1\end{pmatrix} \end{pmatrix}
\end{equation}
where $v_j$ are the first $n-1$ eigenvectors of $\gamma$, chosen as described earlier, 
 maps $(0,...,0,1)$ to $X$, 
maps $(0,...,0,-1)$ to $Y$, and maps $\tilde{C}_0$ to $\tilde{C}$. 
The automorphism 
\begin{equation}
\label{gamma0matrix}
\gamma_0=A_{\gamma}^{-1}\gamma A_{\gamma}=\begin{pmatrix} I_{\gamma} & 0 & 0 \\ 
0 & \dfrac{1}{2\lambda}+\dfrac{\lambda}{2} & -\dfrac{1}{2\lambda}+\dfrac{\lambda}{2} \\ 
0 & -\dfrac{1}{2\lambda}+\dfrac{\lambda}{2} & \dfrac{1}{2\lambda}+\dfrac{\lambda}{2} 
\end{pmatrix}
\end{equation}
leaves $\tilde{C}_0$ invariant.
\begin{remark} 
\label{remJreal}
For each $\xi\in\tilde{C}_0$ $J(A_{\gamma},\xi)$ is real. This follows from 
(\ref{jacobidet}), (\ref{tC0def}), (\ref{Agammamatrix}), and the assumption that $X$ and $Y$ are real. 
\end{remark}
Denote by $C$ the simple closed geodesic $\tilde{C}/\langle \gamma \rangle $ in $M$. We have:  
$\gamma_0:0\mapsto (0,...,0,\dfrac{\lambda^2-1}{\lambda^2+1})$, and  therefore 
the hyperbolic length $l(C)$ of $C$ satisfies 
$$
\cosh \frac{l(C)}{2}=\frac{1+\lambda^2}{2|\lambda|}. 
$$
It follows that that 
\begin{equation}
\label{lnlambda}
\ln |\lambda| = \dfrac{l(C)}{2}. 
\end{equation}

Denote by $K_M$ the canonical bundle on $M$ and denote by $K_{\B^n}$ the canonical bundle on $\B^n$.  
The complex vector space $H^0(M,K_M)$ is isomorphic to the space $H^0_{\Gamma}(\B^n,K_{\B^n})$ of $\Gamma$-invariant holomorphic $n$-forms 
on $\B^n$, i.e. $n$-forms $f(z)dz_1\wedge...\wedge dz_n$, where $f:\B^n\to\C$ is holomorphic and such that  
$$
f(A z) J(A,z)=f(z) \ {\mathrm{for \ all}} \ A\in \Gamma,z\in \B^n.
$$ 
This complex vector space 
has an inner product defined by 
$$
(f(z)dz_1\wedge...\wedge dz_n, g(z)dz_1\wedge...\wedge dz_n)=(-1)^{\frac{n(n-1)}{2}}\Bigl ( \frac{i}{2}\Bigr ) ^n\int \limits_M f(z)\overline{g(z)}
dz_1\wedge...\wedge dz_n\wedge 
d\bar{z}_1\wedge...\wedge d\bar{z}_n=
$$
\begin{equation}
\label{inprodcanb}
\int \limits_M \frac{f(z)\overline{g(z)}}{K(z,z)}dV(z). 
\end{equation}
The complex vector space $H^0(M,K_M^{\otimes k})$ is isomorphic to the space $H^0_{\Gamma}(\B^n,K_{\B^n}^{\otimes k})$ 
of $\Gamma$-invariant holomorphic sections of $K_{\B^n}^{\otimes k}$. Those can be written as $f(z)( dz_1\wedge...\wedge dz_n)^{\otimes k}$, 
 where $f:\B^n\to\C$ is holomorphic and such that  
\begin{equation}
\label{aftransf}
f(A z) J(A,z)^k=f(z) \ {\mathrm{for \ all}} \ A\in \Gamma,z\in \B^n.
\end{equation}
Denote the space of holomorphic functions on $\B^n$ that satisfy (\ref{aftransf}) by  $\tilde{H}^0_{\Gamma}(\B^n,K_{\B^n}^{\otimes k})$. 
This space is isomorphic to $H^0_{\Gamma}(\B^n,K_{\B^n}^{\otimes k})$ and to $H^0(M,K_M^{\otimes k})$. 
The inner product on $\tilde{H}^0_{\Gamma}(\B^n,K_{\B^n}^{\otimes k})$ is  defined by 
$$
(f, g)=
\int \limits_M \frac{f(z)\overline{g(z)}}{K(z,z)^k}dV(z). 
$$
The following statement will be useful. 
\begin{lemma}
\label{lemmainvongeod}
Suppose $f$ is a holomorphic function on $\B^n$ that satisfies (\ref{aftransf}). 
Then the function $f(z)K(z,z)^{-\frac{k}{2}}$, restricted to $\tilde{C}$, is $\gamma$-invariant. 
\end{lemma}
\noindent {\bf Proof.} 
It is sufficient to show that $J(\gamma,z)$ is real-valued for all $z\in \tilde{C}$ (then the statement follows from (\ref{bktransf}) 
and (\ref{aftransf})). For $z\in\tilde{C}$, by (\ref{cocyclecond}) 
$$
J(\gamma,z)=\frac{J(A_{\gamma},\gamma_0w)J(\gamma_0,w)}{J(A_{\gamma},w)}. 
$$
where $w=A_{\gamma}^{-1}z\in \tilde{C}_0$. The statement now follows from Remark \ref{remJreal}, 
(\ref{jacobidet}), and (\ref{gamma0matrix}). $\Box$

For $w\in\B^n$, $k\ge 2$, the series  
\begin{equation}
\label{thetapointdef}
\Theta_w^{(k)}(z)=c(\B^n ,k)
\sum_{A\in\Gamma} \K (A z,w)^kJ(A ,z)^k 
\end{equation}
converges absolutely and uniformly on compact sets,  by \cite[Prop.1 p.44]{baily:73}. We have: 
$\Theta_w^{(k)}\in \tilde{H}^0_{\Gamma}(\B^n,K_{\B^n}^{\otimes k})$, 
and for every $g\in \tilde{H}^0_{\Gamma}(\B^n,K_{\B^n}^{\otimes k})$
\begin{equation}
\label{inprodpoint}
(g, \Theta_w^{(k)})=g(w)
\end{equation}
(this follows from (\ref{bktransf}) and (\ref{wreprod})). 

The $1$-form 
$$
\varphi(z)=\K(z,z)^{\frac{1}{n+1}}dz_n,
$$ 
restricted to $\tilde{C_0}$, is $\gamma_0$-invariant  (this follows from (\ref{jacobidet}), 
(\ref{bktransf}), (\ref{gamma0matrix})).  Then the $1$-form $(A_{\gamma}^{-1})^*\varphi$, restricted to $\tilde{C}$, is $\gamma$-invariant, 
and thus descends to $C$. Since 
$$
\Theta_{hw}^{(k)}( z) \overline{J(h,w)}^k=\Theta_{w}^{(k)}( z) \ {\mathrm{for \ all}} \ h\in \Gamma,w,z\in \B^n,
$$
Lemma \ref{lemmainvongeod} implies that $\Theta_w^{(k)}(z)\K (w,w)^{-\frac{k}{2}}$, regarded as a function of $w$, restricted 
to $\tilde{C}$,  
is $\gamma$-invariant. 

Define the function $\Theta_C^{(k)}$ on $\B^n$ by 
\begin{equation}
\label{thetaCdef}
\Theta_C^{(k)}(z)=\int\limits_C \Theta_w^{(k)}(z)\K (w,w)^{-\frac{k}{2}} ((A_{\gamma}^{-1})^*\varphi) (w).
\end{equation}
It is a holomorphic function, and, moreover, 
$\Theta_C^{(k)}\in \tilde{H}^0_{\Gamma}(\B^n,K_{\B^n}^{\otimes k})$. 
For every $g\in \tilde{H}^0_{\Gamma}(\B^n,K_{\B^n}^{\otimes k})$
\begin{equation}
\label{inprodgeod}
(g, \Theta_C^{(k)})=
\int\limits_C g(z)\K (z,z)^{-\frac{k}{2}} ((A_{\gamma}^{-1})^*\varphi )(z).
\end{equation}
\begin{theorem}
\label{mainth}
$$
( \Theta_{C}^{(k)},\Theta_{C}^{(k)})\sim  k^{n-\frac{1}{2}} \frac{(n+1)^{n-\frac{1}{2}}}
{\pi^{\frac{3n-1}{2n+2}}(n!)^{\frac{n-1}{n+1}}}\frac{l(C)}{\sqrt{2}}
$$
as $k\to\infty$. 
\end{theorem}
\begin{corollary} 
\label{nonvcor}
The function $\Theta_{C}^{(k)}$ is not identically zero for sufficiently large $k$.
\end{corollary}

\subsection{Additional remarks}
\label{sec:addrem}

To put the concepts from subsection \ref{sec:setting} in a more general context, we will now offer some remarks. 

First, there is a standard way of associating a section of a line bundle to a compact Bohr-Sommerfeld Lagrangian submanifold of a compact 
K\"ahler manifold (the Hermitian holomorphic line bundle is related to the K\"ahler form by requiring that 
the curvature of the Chern connection is $-i$ times the K\"ahler form). 
For the $k$-th power of the line bundle, as $k\to\infty$, the square of the norm of this section is asymptotic to 
const$\cdot k^{\frac{n}{2}}$, where $n$ is the complex dimension of the manifold \cite{borthwick:95}. The same procedure can be applied to 
isotropic submanifolds. The closed geodesic $C$  is an isotropic submanifold of $M$, and it is also is a Bohr-Sommerfeld set in $M$ 
(in the terminology of \cite[p.1271]{burns:10}). 
The expression $\dfrac{f(z)\overline{g(z)}}{\K(z,z)}$ in (\ref{inprodcanb}) 
represents the pointwise Hermitian inner product in the holomorphic Hermitian 
line bundle $K_M$. The Chern connection $\nabla$, in a local holomorphic frame $e(z)$, is given by the $1$-form 
$\Theta = \partial \log \frac{e(z)\overline{e(z)}}{\K(z,z)}$. 
Denote by $\iota: C\to M$ the inclusion map.  Now we will define a covariant constant section $\tau$ of $\iota^* K_M^*\Bigr |_C$ 
(existence of such $\tau$ means that a Bohr-Sommerfeld condition is satisfied). 
Let $s=f(z)dz_1\wedge...\wedge dz_n\in H^0_{\Gamma}(\B^n,K_{\B^n})$. 
For $w\in \tilde{C}$ set 
\begin{equation}
\label{taudef}
\tau(s)(w)=f(w)\K(w,w)^{-\frac{1}{2}}.
\end{equation}
Because of Lemma \ref{lemmainvongeod}, 
the equation (\ref{taudef}) defines a section of $\iota^* K_M^*\Bigr |_C$.
To verify that it is covariant constant (i.e. $\nabla^* \tau=0$ over $C$,  where $\nabla^*$ is the connection in $K_M^*$), 
we do this calculation: 
in a local holomorphic frame $e(z)$ for $K_M$, 
for a holomorphic section $s$ of $K_M$ represented locally by $\psi(z)e(z)$ where $\psi$ is a local holomorphic function, for $z\in \tilde{C}$  
$$
(\nabla^* \tau) (s)=d(\tau(s))-\tau(\nabla s)=
$$
$$
d(\psi(z)e(z)K(z,z)^{-\frac{1}{2}})-\K(z,z)^{-\frac{1}{2}}(e(z)d(\psi(z))+\psi(z)e(z)\partial \log \frac{e(z)\overline{e(z)}}{\K(z,z)})=
$$
$$
\K(z,z)^{-\frac{1}{2}}\psi de-
\frac{1}{2}e\psi \K(z,z)^{-\frac{3}{2}}d\K(z,z)
-\K(z,z)^{-\frac{1}{2}} \psi e\partial\log e 
+ \K(z,z)^{-\frac{1}{2}}\psi e\partial\log \K(z,z)=0
$$
since $\partial\log e =d\log e$ (because $e$ is holomorphic) and 
$\partial\log\K(z,z)\Bigr |_{\tilde{C}}=\frac{1}{2}d\log\K(z,z)\Bigr |_{\tilde{C}}$.

We have continuous linear functionals on $H^0(M,K_M^{\otimes k})$: 
$$
s\mapsto \tau ^{\otimes k}(s(w)), {\mbox{for }} \ w\in C
$$
$$
s\mapsto \int\limits _C \tau^{\otimes k}(s(w))  ((A_{\gamma}^{-1})^*\varphi)(w). 
$$
By Riesz representation theorem, there are unique elements  $s_w^{(k)}$, $s_C^{(k)}$ of $H^0(M,K_M^{\otimes k})$ such that 
$$
\tau ^{\otimes k}(s(w))=( s,s_w^{(k)} ) , {\mbox{for }} \ w\in C
$$
$$
\int\limits _C \tau^{\otimes k}(s(w)) ((A_{\gamma}^{-1})^*\varphi)(w)=( s,s_C^{(k)} ). 
$$
Because of (\ref{inprodpoint}), (\ref{inprodgeod}), the element of $\tilde{H}^0_{\Gamma}(\B^n,K_{\B^n}^{\otimes k})$ 
corresponding to  $s_C^{(k)}$ (under the isomorphism $H^0(M,K_M^{\otimes k})\cong 
\tilde{H}^0_{\Gamma}(\B^n,K_{\B^n}^{\otimes k})$)
is $\Theta_C^{(k)}$. 

To summarize, $\Theta_C^{(k)}(z)(dz_1\wedge...\wedge dz_n)^{\otimes k}$ is the $\Gamma$-invariant holomorphic section of 
$K_{\B^n}^{\otimes k}$ that is 
associated to $C$, an isotropic submanifold of $M$ satisfying a Bohr-Sommerfeld condition.

Second, from the perspective of automorphic forms in several complex variables, the function $\Theta_w^{(k)}$ is a Poincar\'e series, 
and the function 
$\Theta_C^{(k)}$ is a relative Poincar\'e series. Indeed, 
$$ 
\Theta_C^{(k)}(z)=c(\B^n ,k)
\sum\limits_{A\in \langle \gamma\rangle\backslash \Gamma}\int\limits_{\tilde{C}}\K (A z,w)^k \K (w,w)^{-\frac{k}{2}} 
((A_{\gamma}^{-1})^*\varphi)(w)J(A,z)^k.
$$ 
Relative Poincar\'e series associated to closed geodesics in ball quotients have been previously studied in \cite{katok:85, katok:87} and \cite{borthwick:95}
for $n=1$ (for compact Riemann surfaces of genus $g\ge 1$), 
and for $n\ge 1$ in the publications \cite{foth:01, foth:04} co-authored by S. Katok and T. Foth (T. Barron). 
A somewhat tedious calculation shows that $\Theta_C^{(k)}$ is equal, up to a nonzero constant depending on $n$, $k$, $X$, $Y$, to the function 
${\mathcal{P}}_C^{(k)}$ defined by 
$$
{\mathcal{P}}_C^{(k)}(z)=\sum\limits_{A\in \langle \gamma\rangle\backslash \Gamma} \dfrac{1}{(\langle A z,X\rangle \langle A z,Y\rangle )^
{\frac{(n+1)k}{2}} }J(A,z)^{k}
$$
when $(n+1)k$ is an even integer (so, $n$ is odd or $k$ is even). The relative Poincar\'e series ${\mathcal{P}}_C^{(k)}$ were used in  
\cite{foth:01} and \cite{foth:04} to address the spanning question for the space of cusp forms. 
We study a different question, about the asymptotics of the inner products. The normalizing factor 
in the definition of  $\Theta_C^{(k)}$ is determined by (\ref{inprodgeod}). 
The statement about the $k\to\infty$ asymptotics of $(\Theta_C^{(k)},\Theta_C^{(k)})$ is Theorem \ref{mainth}. One important consequence is a nonvanishing result, 
Corollary \ref{nonvcor}. (In general the Poincar\'e series map 
has a large kernel, which leads to the question about nonvanishing of Poincar\'e series. A similar question 
can be posed for relative Poincar\'e series.)

\subsection{Proof of Theorem \ref{mainth}} From (\ref{inprodgeod}) 
$$
( \Theta_{C}^{(k)},\Theta_{C}^{(k)}) =\int\limits _{C} \Theta_{C}^{(k)}(z)\K (z,z)^{-\frac{k}{2}}((A_{\gamma}^{-1})^*\varphi )(z)=
\int\limits _0^{\gamma_0(0)} \Theta_{C}^{(k)}(A_{\gamma}w) \K (A_{\gamma}w,A_{\gamma}w)^{-\frac{k}{2}} \ \varphi(w)
$$
where the line integral $\int\limits _0^{\gamma_0(0)}$ is over the segment of $\tilde{C}_0$ that connects $0$ and $\gamma_0(0)$. 
Using (\ref{cocyclecond}), (\ref{bktransf}), (\ref{thetapointdef}), (\ref{thetaCdef}), and Remark \ref{remJreal}, we get: 
$$
( \Theta_{C}^{(k)},\Theta_{C}^{(k)}) =\int\limits _0^{\gamma_0(0)}   \int\limits _{C} \Theta_{\zeta}^{(k)}(A_{\gamma}w)\K (\zeta,\zeta)^{-\frac{k}{2}}((A_{\gamma}^{-1})^*\varphi)(\zeta)                       
\K(A_{\gamma}w,A_{\gamma}w)^{-\frac{k}{2}} \ \varphi(w)=
$$
$$
\int\limits _0^{\gamma_0(0)}\int\limits _0^{\gamma_0(0)}
\Theta_{A_{\gamma}\xi}^{(k)}(A_{\gamma}w)\K (A_{\gamma}\xi,A_{\gamma}\xi)^{-\frac{k}{2}} \ \varphi (\xi)                       
\K(A_{\gamma}w,A_{\gamma}w)^{-\frac{k}{2}} \ \varphi(w)=
$$
$$
c(\B^n ,k)\int\limits _0^{\gamma_0(0)}\int\limits _0^{\gamma_0(0)}
\sum_{g\in\Gamma} \K (g A_{\gamma}w,A_{\gamma}\xi)^k J(g,A_{\gamma}w)^k 
\K (A_{\gamma}\xi ,A_{\gamma}\xi)^{-\frac{k}{2}} \ \varphi(\xi) 
\K (A_{\gamma}w,A_{\gamma}w)^{-\frac{k}{2}}  \ \varphi(w)=
$$
$$
c(\B^n ,k)\int\limits _0^{\gamma_0(0)}\int\limits _0^{\gamma_0(0)}
\sum_{h\in A_{\gamma}^{-1}\Gamma A_{\gamma}} \K (hw,\xi)^k J(h,w)^k \K (\xi ,\xi)^{-\frac{k}{2}} \ \varphi(\xi) 
\K (w,w)^{-\frac{k}{2}} \ \varphi(w).
$$
Denote
$$
J_1(k)=c(\B^n ,k)\int\limits _0^{\gamma_0(0)}\int\limits _0^{\gamma_0(0)}
\sum_{\substack{{h\in A_{\gamma}^{-1}\Gamma A_{\gamma}}\\{h\notin \langle \gamma_0\rangle}}} 
\K (hw,\xi)^k J(h,w)^k \K (\xi ,\xi)^{-\frac{k}{2}} \ \varphi(\xi) 
\K (w,w)^{-\frac{k}{2}} \ \varphi(w)
$$
We will show that 
$|J_1(k)|$  
is rapidly decreasing as $k\to\infty$. We have: 
$$
|J_1(k)|=
|c(\B^n ,k)\int\limits _0^{\gamma_0(0)}\int\limits _0^{\gamma_0(0)}
\sum_{\substack{{h\in A_{\gamma}^{-1}\Gamma A_{\gamma}}\\{h\notin\langle \gamma_0\rangle}}} \K (h w,\xi)^k \K (\xi,\xi)^{-\frac{k}{2}} 
 \ \varphi(\xi)J(h ,w) ^k 
\K (w,w)^{-\frac{k}{2}} \ \varphi(w)|\le 
$$
$$
c(\B^n ,k)\int\limits_{0}^{\frac{\lambda^2-1}{\lambda^2+1}} \int\limits _{0}^{\frac{\lambda^2-1}{\lambda^2+1}}
\sum_{\substack{{h\in A_{\gamma}^{-1}\Gamma A_{\gamma}}\\{h\notin\langle\gamma_0\rangle}}} |\K (hw,\xi)|^k \K (\xi,\xi)^{-\frac{k}{2}+\frac{1}{n+1}} \ d({\mathrm{Re}} \ \xi_n)|J(h ,w)|^k 
\K (w,w)^{-\frac{k}{2}+\frac{1}{n+1}} \ d({\mathrm{Re}} \ w_n).
$$
Since $\Gamma$ is discrete, there is $\delta_0>0$ such that the hyperbolic distance between $hw$ and $\xi$  
is not less than $\delta_0$ (for all $w\in\tilde{C}_0$ on the segment between $0$ and $\gamma_0(0)$ and all 
$h\in A_{\gamma}^{-1}\Gamma A_{\gamma}$ such that $h\notin\langle\gamma_0\rangle$). Hence 
$$
\cosh ^2 \frac{\rho(hw,\xi)}{2}=\frac{\langle hw,\xi\rangle \langle \xi,hw\rangle}{\langle hw,hw\rangle \langle \xi,\xi\rangle}\ge \cosh^2\frac{\delta_0}{2}.
$$
Then in the integral above, for $k\ge 3$, 
$$
|\K (hw,\xi)|^k= \Bigl ( \frac{n!}{\pi^n}\Bigr )^k \frac{1}{|\langle hw,\xi\rangle |^{(n+1)k}}\le 
\Bigl ( \frac{n!}{\pi^n}\Bigr )^k \frac{1}{|\langle hw,\xi\rangle |^{2(n+1)}\Bigl [ 
\langle hw,hw\rangle \langle \xi,\xi\rangle \cosh^2\frac{\delta_0}{2} \Bigr ]^{\frac{(n+1)k}{2}-(n+1)} }=
$$
$$
\frac{|\K(hw,\xi)|^2 \K(hw,hw)^{\frac{k}{2}-1} \K(\xi,\xi)^{\frac{k}{2}-1}}{(\cosh^2\frac{\delta_0}{2})^{\frac{(n+1)k}{2}-(n+1)}}=
 \frac{|\K(hw,\xi)|^2 \K(w,w)^{\frac{k}{2}-1}\K(\xi,\xi)^{\frac{k}{2}-1}}{(\cosh^2\frac{\delta_0}{2})^{\frac{(n+1)k}{2}-(n+1)}
 |J(h,w)|^{k-2}}.
$$
We get: 
$$
|J_1(k)|\le 
c(\B^n ,k)\frac{1}{(\cosh^2\frac{\delta_0}{2})^{\frac{(n+1)k}{2}-(n+1)}}
\int\limits _{0}^{\frac{\lambda^2-1}{\lambda^2+1}} \int\limits _{0}^{\frac{\lambda^2-1}{\lambda^2+1}}
\sum_{\substack{{h\in A_{\gamma}^{-1}\Gamma A_{\gamma}}\\{h\notin\langle \gamma_0\rangle}}} |\K (hw,\xi)|^2| \K (\xi,\xi)^{-\frac{n}{n+1}}
d({\mathrm{Re}} \ \xi_n) 
$$
$$
|J(h,w)|^2
\K (w,w)^{-\frac{n}{n+1}} d({\mathrm{Re}} \ w_n). 
$$
The series $\sum\limits_{h\in A_{\gamma}^{-1}\Gamma A_{\gamma}} |\K (hw,\xi)|^2 |J(h,w)|^2$ (with a fixed $\xi$) converges uniformly 
on compact subsets of $\B^n$ by \cite[Prop.1 p.44]{baily:73}. For $\xi\in \tilde{C}_0$ $K(\xi,\xi)^{-\frac{n}{n+1}}=
\Bigl ( \dfrac{n!}{\pi^n}\Bigr )^{-\frac{n}{n+1}}(1-({\mathrm{Re}} \ \xi_n)^2)^n$. 
We conclude that 
$$
|J_1(k)|\le c(\B^n ,k)\frac{{\mathrm{const}}(\lambda,n)}{(\cosh^2\frac{\delta_0}{2})^{\frac{(n+1)k}{2}-(n+1)}}.
$$
Since $c(\B^n ,k)\sim {\mathrm{const}}(n) k^n$ and $\cosh^2\dfrac{\delta_0}{2}>1$, 
this implies that for any $l\in N$ there is $C>0$ such that $|J_1(k)|<\dfrac{C}{k^l}$ as $k\to\infty$. 

It remains to estimate 
$$
J_2(k)=c(\B^n ,k)\int\limits _0^{\gamma_0(0)}\int\limits _0^{\gamma_0(0)}
\sum_{h\in \langle \gamma_0\rangle} 
\K (hw,\xi)^k  \K (\xi ,\xi)^{-\frac{k}{2}} J(h,w)^k
\K (w,w)^{-\frac{k}{2}} \varphi (\xi) \varphi (w).
$$
By (\ref{bktransf}) and Remark \ref{remJreal}, in the integral above, $\K (w ,w)^{-\frac{k}{2}} J(h,w)^k=\K (hw,hw)^{-\frac{k}{2}}$. 
Denote  $u_n={\mathrm{Re}} \ \xi _n$ and denote $z=hw$,  $x_n={\mathrm{Re}} \  z_n$. 
$$
J_2(k)=c(\B^n ,k)\int\limits_0^{\frac{\lambda^2-1}{\lambda^2+1}}\int\limits_{-1}^1
\K (z,\xi)^k \K (\xi,\xi)^{-\frac{k}{2}+\frac{1}{n+1}} 
\K (z,z)^{-\frac{k}{2}+\frac{1}{n+1}} \ dx _n \ du_n =
$$
$$
c(\B^n ,k)\Bigl ( \frac{n!}{\pi^n} \Bigr ) ^{\frac{2}{n+1}}
\int\limits_{0}^{\frac{\lambda^2-1}{\lambda^2+1}}
\int\limits_{-1}^{1}
\frac{(\langle \xi,\xi\rangle \langle z,z\rangle )^{\frac{(n+1)k}{2}-1}}{(-\langle z,\xi\rangle)^{(n+1)k}}  \ dx_n 
\ du_n =
$$
$$
c(\B^n ,k)\Bigl ( \frac{n!}{\pi^n} \Bigr ) ^{\frac{2}{n+1}}
\int\limits_{0}^{\frac{\lambda^2-1}{\lambda^2+1}}
\int\limits_{-1}^{1}\frac{(1-x_n^2)^{\frac{(n+1)k}{2}-1}}{(1-x_nu_n)^{(n+1)k}}dx_n(1-u_n^2)^{\frac{(n+1)k}{2}-1}du_n.
$$
We will apply the Laplace method \cite[II.1,(1.5)]{wong:89}, \cite{hsu:56}  to the integral  
$\int\limits_{-1}^{1}\dfrac{(1-x_n^2)^{\frac{(n+1)k}{2}-1}}{(1-x_nu_n)^{(n+1)k}}dx_n$, with a fixed $u_n$. 
$$
\int\limits_{-1}^{1}\frac{(1-x_n^2)^{\frac{(n+1)k}{2}-1}}{(1-x_nu_n)^{(n+1)k}}dx_n=
\int\limits_{-1}^1f(x_n)^{(n+1)k-2}\frac{1}{(1-x_nu_n)^2}dx_n
$$
where $f(x_n)=\frac{\sqrt{1-x_n^2}}{1-x_nu_n}$. We have: 
$$
f'(x_n)=\frac{u_n-x_n}{\sqrt{1-x_n^2}(1-x_nu_n)^2}, \ f''(x_n)\Bigr |_{x_n=u_n}=-\frac{1}{(1-u_n^2)^{\frac{5}{2}}}<0,
$$
$$
\int\limits_{-1}^1f(x_n)^{(n+1)k-2}\frac{1}{(1-x_nu_n)^2}dx_n\sim \frac{1}{(1-u_n^2)^2}f(u_n)^{(n+1)k-2+\frac{1}{2}}\Bigl ( 
\frac{-2\pi}{((n+1)k-2)f''(u_n) }   \Bigr ) ^{\frac{1}{2}}=
$$
$$
\sqrt{\frac{2\pi}{(n+1)k-2}}(1-u_n^2)^{-\frac{(n+1)k}{2}}.
$$
Therefore, as $k\to\infty$, 
$$
J_2(k)\sim c(\B^n ,k) \Bigl ( \frac{n!}{\pi^n} \Bigr ) ^{\frac{2}{n+1}}  \frac{\sqrt{2\pi}}{\sqrt{(n+1)k-2}}
\int\limits_0^{\frac{\lambda^2-1}{\lambda^2+1}}
(1-u_n^2)^{-1}du_n=
$$
$$
c(\B^n ,k) \Bigl ( \frac{n!}{\pi^n} \Bigr ) ^{\frac{2}{n+1}}  \frac{\sqrt{2\pi}}{\sqrt{(n+1)k-2}}\ln|\lambda| .
$$
Then, using  (\ref{ballconstk}), we get: 
$$
J_2(k)
\sim
k^{n-\frac{1}{2}} \frac{\sqrt{2} \ (n+1)^{n-\frac{1}{2}}}
{\pi^{\frac{3n-1}{2n+2}}(n!)^{\frac{n-1}{n+1}}}\ln|\lambda|
$$
as $k\to\infty$. The statement now follows from (\ref{lnlambda}). 
$\Box$
\begin{remark}
The next term in the asymptotic expansion of $( \Theta_{C}^{(k)},\Theta_{C}^{(k)})$ is of the form 
${\mathrm{const}}\cdot k^{n-\frac{3}{2}}$. This follows from the proof of Theorem \ref{mainth}. Indeed,   
$( \Theta_{C}^{(k)},\Theta_{C}^{(k)})=J_1(k)+J_2(k)$, $|J_1(k)|$ is rapidly decreasing as $k\to\infty$,
and $J_2(k)$ is equal to $\Bigl ( \dfrac{n!}{\pi^n} \Bigr ) ^{\frac{2}{n+1}}$, times $c(\B^n ,k)$, times an integral. 
Multiplying (\ref{ballconstk}) and the asymptotic expansion   
for the integral (taking into account the first two terms from the Laplace method), we obtain the conclusion.   
\end{remark}
\begin{remark}
Theorem \ref{mainth} does not follow from the theorems  in  \cite{alluhaibi:18} related to submanifolds of the ball. 
\end{remark}


\begin{thebibliography}{999}

\bibitem[AB]{alluhaibi:18}N. Alluhaibi, T. Barron. 
\newblock {\em On vector-valued automorphic forms on bounded symmetric domains.} 
\newblock Preprint, 2018. \  https://arxiv.org/abs/1806.03779   

\bibitem[B]{baily:73}W. Baily.  
\newblock {\em Introductory lectures on automorphic forms.} 
\newblock Iwanami Shoten, Publishers, Tokyo; Princeton University Press, Princeton, N.J., 1973. 

\bibitem[BPU]{borthwick:95}D. Borthwick, T. Paul, A. Uribe.  
\newblock {\em Legendrian distributions with applications to relative Poincar\'e series.} 
\newblock Invent. Math. 122 (1995), no. 2, 359-402.

\bibitem[BGW]{burns:10}D. Burns, V. Guillemin, Z. Wang. 
\newblock {\em Stability functions.} 
\newblock Geom. Funct. Anal. 19 (2010), no. 5, 1258-1295. 

\bibitem[dB]{debru:61}N. de Bruijn. 
\newblock {\em Asymptotic methods in analysis.}
\newblock Second edition. Bibliotheca Mathematica, Vol. IV. 
\newblock North-Holland Publishing Co., Amsterdam; P. Noordhoff Ltd., Groningen 1961. 

\bibitem[DP]{deber:06}M. Debernardi, R. Paoletti. 
\newblock {\em Equivariant asymptotics for Bohr-Sommerfeld Lagrangian submanifolds.}
\newblock Comm. Math. Phys. 267 (2006), no. 1, 227-263. 

\bibitem[F1]{foth:02}T. Foth.
\newblock {\em Bohr-Sommerfeld tori and relative Poincar\'e series on a complex hyperbolic space.}
\newblock Commun. Anal. Geom. 10, no. 1, 151 (2002). 

\bibitem[F2]{foth:08}T. Foth. 
\newblock {\em Legendrian tori and the semi-classical limit.} 
\newblock Diff. Geom. Appl.  26, no. 1, 63 (2008).

\bibitem[FK1]{foth:01}T. Foth, S. Katok.
\newblock {\em Spanning sets for automorphic forms and dynamics of the frame flow on complex hyperbolic spaces.}
\newblock Ergodic Theory Dynam. Systems 21 (2001), no. 4, 1071-1099.

\bibitem[FK2]{foth:04}T. Foth, S. Katok.
\newblock Appendix to {\em S. Katok.
Livshitz theorem for the unitary frame flow.} 
\newblock Ergodic Theory Dynam. Systems 24 (2004), no. 1, 127-140; pp. 137-140. 

\bibitem[G]{goldman:99}W. Goldman. 
\newblock {\em Complex hyperbolic geometry.}
\newblock Oxford Mathematical Monographs. Oxford Science Publications. 
\newblock The Clarendon Press, Oxford University Press, New York, 1999.

\bibitem[GUW]{guillemin:16}V. Guillemin, A. Uribe, Z. Wang. 
\newblock  {\em Semiclassical states associated with isotropic submanifolds of phase space.} 
\newblock    Lett. Math. Phys.  106 (2016), no. 12, 1695-1728. 

\bibitem[H]{hsu:56}L. Hsu. 
\newblock {\em On an asymptotic integral.}
\newblock  Proc. Edinburgh Math. Soc. (2) 10 (1956), 141-144. 

\bibitem[JW]{jeffrey:92}L. Jeffrey, J. Weitsman. 
\newblock {\em Bohr-Sommerfeld orbits in the moduli space of flat connections and the Verlinde dimension formula.} 
\newblock Comm. Math. Phys. 150 (1992), no. 3, 593-630. 

\bibitem[Ka]{katok:85}S. Katok. 
\newblock {\em Closed geodesics, periods and arithmetic of modular forms.}
\newblock  Invent. Math. 80 (1985), 469-480.

\bibitem[KaM]{katok:87}S. Katok, J. Millson. 
\newblock {\em Eichler-Shimura homology, intersection numbers and rational structures on spaces of modular forms.} 
\newblock Trans. Amer. Math. Soc. 300 (1987), no. 2, 737-757. 

\bibitem[Kr]{krantz:01}S. Krantz. 
\newblock {\em Function theory of several complex variables.}
Reprint of the 1992 edition. AMS Chelsea Publishing, Providence, RI, 2001. 

\bibitem[KuM]{kudla:79}S. Kudla, J. Millson. 
\newblock {\em Harmonic differentials and closed geodesics on a Riemann surface.}
\newblock Invent. Math. 54 (1979), no. 3, 193-211. 

\bibitem[P]{paoletti:08}R. Paoletti. 
\newblock {\em A note on scaling asymptotics for Bohr-Sommerfeld Lagrangian submanifolds.}
\newblock Proc. Amer. Math. Soc. 136 (2008), no. 11, 4011-4017. 

\bibitem[R]{rudin:80}W. Rudin. 
\newblock {\em Function theory in the unit ball of $\C^n$.}
\newblock  Springer-Verlag, New York-Berlin, 1980.

\bibitem[S]{selberg:57}A. Selberg. 
\newblock {\em Automorphic functions and integral operators.} 
\newblock In {\em Collected papers}, vol. I,  Springer-Verlag, 1989; 464-468.

\bibitem[TW]{tong:83}Y. Tong, S. Wang. 
\newblock {\em Theta functions defined by geodesic cycles in quotients of SU(p,1).} 
\newblock Invent. Math. 71 (1983), no. 3, 467-499. 

\bibitem[W]{wong:89}R. Wong. 
\newblock {\em  Asymptotic approximations of integrals.} 
\newblock  Academic Press, Inc., Boston, MA, 1989. 

\bibitem[Z]{zhu:05}K. Zhu. 
\newblock {\em Spaces of holomorphic functions in the unit ball.} 
\newblock Springer-Verlag, New York, 2005. 

\end{thebibliography}
\end{document}